\documentclass[12pt]{article}

\newtheorem{theorem}{Theorem}

\newtheorem{corollary}[theorem]{Corollary}

\newtheorem{example}[theorem]{Example}

\newtheorem{remark}[theorem]{Remark}

\newenvironment{proof}[1][Proof]{\textbf{#1.} }{\ \rule{0.5em}{0.5em}}
\newdimen\dummy
\dummy=\oddsidemargin
\date{}
\begin{document}

\title{A retraction theorem for topological fundamental groups with application to
the Hawaiian earring}
\author{Paul Fabel \\
Department of Mathematics \& Statistics\\
Mississippi State University}
\maketitle

\begin{abstract}
A characterization of regular topological fundamental groups yields a `no
retraction theorem' for spaces constructed in similar fashion to the
Hawaiian earring.
\end{abstract}

\section{Introduction}

Analysis of algebraic properties of fundamental groups $\pi _{1}(X)$ and $%
\pi _{1}(Y)$ often leads to the conclusion that $X$ cannot be embedded as a
retract of $Y.$ This paper illustrates the possibility to reach the 
same conclusion by analyzing \textit{topological}\textbf{\ }properties of 
these groups. 

Theorem \ref{main} offers a characterization of those spaces $Y,$ whose
fundamental group $\pi _{1}(Y),$ topologized in a natural way, is a $T_{1}$
space: The space $\pi _{1}(Y)$ is $T_{1}$ if and only if each retraction of $%
Y$ induces an embedding between topological fundamental groups.

Recent work of the author \cite{fab2} shows that the topological fundamental
group of a space $X$ constructed in similar fashion to the Hawaiian earring
is not a Baire space. 

Consequently (Theorem \ref{app}), $X$ cannot be embedded as a retract of a
space $Y$ whose topological fundamental group is completely metrizable. For
example $X$ cannot be emdedded as a retract of the countable product of
locally simply connected spaces.

\section{Definitions and preliminaries}

Suppose $X$ is a metrizable space and $p\in X.$ Let $C_{p}(X)=\{f:[0,1]%
\rightarrow X$ such that $f$ is continuous and $f(0)=f(1)=p\}.$ Endow $%
C_{p}(X)$ with the topology of uniform convergence.

The \textbf{topological fundamental group} $\pi _{1}(X,p)$ is the set of
path components of $C_{p}(X)$ endowed with the quotient topology under the
canonical surjection $q:C_{p}(X)\rightarrow \pi _{1}(X,p)$ satisfying $%
q(f)=q(g)$ if and only of $f$ and $g$ belong to the same path component of $%
C_{p}(X).$

Thus a set $U\subset \pi _{1}(X,p)$ is open in $\pi _{1}(X,p)$ if and only
if $q^{-1}(U)$ is open in $\pi _{1}(X,p).$

\begin{remark}
The space $\pi _{1}(X,p)$ is a topological group under concatenation of
paths. (Proposition 3.1\cite{Biss}). A map $f:X\rightarrow Y$ determines a
continuous homomorphism $f^{\ast }:\pi _{1}(X,p)\rightarrow \pi _{1}(Y,f(p))$
via $f^{\ast }([\alpha ])=[f(\alpha )]$ (Proposition 3.3 \cite{Biss}).
\end{remark}

Let $[p]\in \pi _{1}(X,p)$ denote the trivial element. Thus $[p]$ is the
path component of the constant path in $C_{p}(X).$

If $X\subset Y$ then $X$ is a \textbf{retract} of $Y$ if there exists a map $%
f:Y\rightarrow X$ such that $f_{X}=id_{X}.$ The space $X$ is $T_{1}$ if each
one point subset of $X$ is closed. The $T_{1}$ space $X$ is \textbf{%
completely regular} if for each closed set $A$ and each point $x^{\ast
}\notin A$ there exists a map $f:X\rightarrow \lbrack 0,1]$ such that $%
f(A)=0 $ and $f(x^{\ast })=1.$ The topological space $X$ is \textbf{%
completely metrizable} if $X$ admits a complete metric compatible with its
topology.

A metric space $(X,d)$ is \textbf{similar to the Hawaiian earring }provided
all of the following hold. Suppose $p\in X$ and $X=\cup _{n=1}^{\infty
}Y_{n} $ with $Y_{n}\cap Y_{m}=\{p\}$ whenever $n\neq m.$ Assume the space $%
Y_{n}$ is path connected, locally simply connected at $p$, and $Y_{n}$ is
not simply connected. Assume also that $\lim_{n\rightarrow \infty
}diam(Y_{n})=0. $ For example if $Y_{n}$ is a simple closed curve then $X$
is the familiar Hawaiian earring.

\section{A characterization of $T_{1}$ topological fundamental groups}

The metrizable spaces with $T_{1}$ topological fundamental groups are
precisely the spaces whose retracts induce embeddings between fundamental
groups.

\begin{theorem}
\label{main}Suppose $Y$ is a metrizable space. The following are equivalent.
\end{theorem}

\begin{enumerate}
\item  The trivial element of $\pi _{1}(Y,p)$ is closed in $\pi _{1}(Y,p).$

\item  $\pi _{1}(Y,p)$ is a $T_{1}$ space.

\item  $\pi _{1}(Y,p)$ is completely regular.

\item  Whenever $X$ is a retract of $Y$ and $j^{\ast }:\pi
_{1}(X,p)\hookrightarrow \pi _{1}(Y,p)$ is the monomorphism induced by
inclusion, then $j^{\ast }$ is an embedding onto a closed subgroup of $\pi
_{1}(Y,p).$
\end{enumerate}

\begin{proof}
The equivalence of $1$ $2$ and $3$ follows from elementary facts about
topological groups (ex. 6 p. 145, ex. 5 p237\cite{Munkres}). 

$1\Rightarrow 4.$ Suppose $r:Y\rightarrow X$ is a retraction. Let $%
q:C_{p}(X)\rightarrow \pi _{1}(X,p)$ and $Q:C_{p}(Y)\rightarrow \pi _{1}(Y,p)
$ denote the canonical quotient maps. Suppose $A\subset \pi _{1}(X,p)$ is
closed. Let $B=\phi (A)\subset \pi _{1}(Y,p).$ To prove $B$ is closed it
suffices to prove $Q^{-1}(B)$ is closed in $C_{p}(Y).$ Suppose $g\in 
\overline{Q^{-1}(B)}.$ Let $g=\lim g_{n}$ with $Q(g_{n})\in B.$ Note $%
r(g_{n})\rightarrow r(g).$ Since $Q(g_{n})\in B$ there exists $f_{n}\in
C_{p}(X)$ path homotopic in $Y$ to $g_{n}$ such that $f_{n}\in q^{-1}(A).$
Thus $r(g_{n})$ and $r(f_{n})=f_{n}$ are path homotopic in $X.$ Thus $%
r(g_{n})\in q^{-1}(A).$ Since $X$ is closed in $Y,$ $C_{p}(X)$ is closed in $%
C_{p}(Y).$ Since $q^{-1}(A)$ is closed in $C_{p}(X),$ and since $C_{p}(X)$
is closed in $C_{p}(Y)$ it must be that $r(g)\in q^{-1}(A).$ Note $g_{n}\ast
r(g_{n})$ is homotopically trivial and converges to $g\ast r(g).$ Since $[p]$
is closed in $\pi _{1}(Y,p)$ it follows that the path component of the
constant map is a closed subspace of $C_{p}(Y).$ Thus $g\ast r(g)$ must be
homotopically trivial in $Y.$ Hence $g$ and $r(g)$ are path homotopic in $Y.$
Thus $g\in Q^{-1}(B).$ Hence $B$ is closed. Therefore $\phi $ is a closed
map and hence an embedding.

$4\Rightarrow 1.$ Note the one point space $X=\{p\}$ is a retract of $Y.$
Thus $j^{\ast }(\pi _{1}(X,p))$ is a closed subspace of $\pi _{1}(Y,p).$
\end{proof}

\begin{remark}
Theorem \ref{main} does not apply to all spaces $Y$. The harmonic
archipelago, explored in detail in \cite{fab2}, provides an example of a
compact path connected metric space $Y$ such that $\pi _{1}(Y,p)$ is not a $%
T_{1}$ space.
\end{remark}

\section{Application: A ``no retraction'' theorem for the Hawaiian earring}

\begin{corollary}
\label{ezcor}Suppose each of $Z$ and $Y$ are metrizable. Suppose $\pi
_{1}(Z,p)$ is not completely metrizable, and suppose $\pi _{1}(Y,p)$ is
completely metrizable. Then $Z$ cannot be embedded as a retract of $Y.$
\end{corollary}

\begin{proof}
Since $\pi _{1}(Y,p)$ is metrizable, the one point subsets of $\pi _{1}(Y,p)$
are closed subspaces. In particular if $X$ is a retract of $Y$ then by
Theorem \ref{main} $\pi _{1}(X,p)$ is homeomorphic to a closed subspace of $%
\pi _{1}(Y,p).$ Thus $\pi _{1}(X,p)$ completely metrizable. Hence $Z$ cannot
be a retract of $Y.$
\end{proof}

\begin{theorem}
\label{app}Suppose $X$ is similar to the Hawaiian earring, suppose $Y$ is
metrizable and suppose $\pi _{1}(Y,p)$ is completely metrizable. Then $X$
cannot be embedded as a retract of $Y.$
\end{theorem}

\begin{proof}
The main result of \cite{fab2} is that $\pi _{1}(X,p)$ is not a Baire space.
Hence $\pi _{1}(X,p)$ is not completely metrizable.
\end{proof}

\begin{example}
Suppose $Y=Z_{1}\times Z_{2}\times ...$where each $Z_{n}$ has the homotopy
type of a bouquet of $n$ loops. Then $Y$ shares some properties with the
Hawaiian earring. For example $Y$ is not locally contractible, for each $%
Z_{n}$ is a retract of $Y,$ and $\pi _{1}(Y)$ is uncountable. However,
(Proposition 5.2 \cite{Biss}) $\pi _{1}(Y)$ is canonically isomorphic and
homeomorphic to the product $\pi _{1}(Z_{1})\times \pi _{1}(Z_{2})...$and
hence completely metrizable (since each factor has the discrete topology.)
Thus $Y$ has no retract similar to the Hawaiian earring.
\end{example}

\end{document}